\documentclass[preprint,times,1p]{elsarticle}
\begin{document}
\begin{frontmatter}
\author{Hal Finkel}
\ead{hal.finkel@yale.edu}
\address{Department of Physics, Yale University, P.O. Box 208120, New Haven, CT 06520-8120}
\title{The differential transformation method and Miller's recurrence}

\begin{abstract}
The differential transformation method (DTM) enables the easy construction of a power-series solution to a nonlinear differential equation. The exponentiation operation has not been specifically addressed in the DTM literature, and constructing it iteratively is suboptimal. The recurrence for exponentiating a power series by J.C.P. Miller provides a concise implementation of exponentiation by a positive integer for DTM. An equally-concise implementation of the exponential function is also provided.
\end{abstract}

\end{frontmatter}

\section{Introduction}
Constructing power-series solutions to differential equations, especially those which do not admit a closed-form solution, has long been an important, and widely-used, solution technique. Traditionally, computing power-series solutions required a fair amount of ``boiler-plate'' symbolic manipulation, especially in the setup of the power-matching phase. The differential transformation method (DTM) enables the easy construction of a power-series solution by specifying a conversion between the differential equation and a recurrence relation for the power-series coefficients~\cite{Chen1996}.

The table in the current literature which specifies the translation between the terms of the differential equation and the recurrence relation has a striking omission: it contains no exponentiation operation. Exponentiation by a positive integer can be constructed iteratively using the table entry for multiplication (i.e. multiplying the function with itself $n$ times), but such a construction is suboptimal because it leads to $n-1$ nested sums. Using a recurrence for exponentiating a power series by J.C.P. Miller, a table entry for positive-integer exponentiation can be provided which introduces only a single sum. A single-sum recurrence for the exponential function of a power series can be similarly constructed.

\section{Preliminaries}
\subsection{Differential transformation}
The differential transform of the function $w(x)$, called $W(k)$, is defined as~\cite{Chen1996}
\begin{equation}
	\label{eq:difftrans}
	W(k) = \frac{1}{k!} \left [ \frac{\partial^k}{\partial x^k} w(x) \right ]_{x=0}
\end{equation}
where $\frac{\partial^k}{\partial x^k}$ is the $k^{\hbox{th}}$ derivative with respect to $x$. The inverse transformation is
\begin{equation}
	\label{eq:invdifftrans}
	w(x) = \sum_{k=0}^{\infty} W(k) x^k \, .
\end{equation}
Combining Equations~\ref{eq:difftrans} and~\ref{eq:invdifftrans} it is clear that the differential transform is derived from the Taylor-series representation of $w(x)$ about $x=0$
\begin{equation}
	w(x) = \sum_{k=0}^{\infty} \frac{1}{k!} \left [ \frac{\partial^k}{\partial x^k} w(x) \right ]_{x=0} x^k \, .
\end{equation}
Table~\ref{tab:basedtm} shows the traditional list of operations and their differential transforms\footnote{Recent literature has also included tables with entries for additional elementary functions (e.g. $\sin(x), \cos(x), e^x$) and integral relations~\cite{Hassan2009}}.

\begin{table}[h]
\begin{center}
\begin{tabular}{|l|l|}
\hline
\bf{Original Function}                        & \bf{Transformed Function}         \\
\hline
$w(x) = y(x) \pm z(x)$                        & $W(k) = Y(k) \pm Z(k)$            \\
$w(x) = \lambda y(x)$, $\lambda$ a constant   & $W(k) = \lambda Y(k)$             \\
$w(x) = \frac{\partial^m}{\partial x^m} y(x)$ & $W(k) = \frac{(k+m)!}{k!} Y(k+m)$ \\
$w(x) = y(x) z(x)$                            & $W(k) = \sum_{l=0}^k Y(l) Z(k-l)$ \\
$w(x) = x^m$                                  & $W(k) = \delta_{k,m}$, $\delta$ is the Kronecker delta \\
\hline
\end{tabular}
\end{center}
\caption{Traditional operations under DTM}
\label{tab:basedtm}
\end{table}

\subsection{J.C.P. Miller's recurrence}
In Knuth's classic book on seminumerical algorithms~\cite{Knuth1981}, he provides a recurrence relation attributed to J.C.P. Miller for the coefficients of a power series raised to some integer power. Unfortunately, the recurrence is not well known~\cite{Zeilberger1995}. Given some formal power series
\begin{equation}
	w(x) = \sum_{k=0}^N a_k x^k
\end{equation}
where $N$ could be $\infty$ and $a_0 \neq 0$, then
\begin{equation}
	w(x)^m = \left ( \sum_{k=0}^N a_k x^k \right )^m = \sum_{k=0}^{Nm} c_k x^k
\end{equation}
where the $\{c_k\}$ are given by the recurrence relation~\cite{Zeilberger1995}
\begin{equation}
	\begin{array}{l}
	c_0 = a_0^m \\
	c_k = \frac{1}{k a_0} \sum_{j=1}^N [(m+1)j - k] a_j c_{k-j} \, .
	\end{array}
\end{equation}
A simple proof by Zeilberger is as follows~\cite{Zeilberger1995}: $c_k$ is the coefficient of $x^0$ in the Laurent series expansion of $w(x)^m/x^k$. For any Laurent series, $f(x)$, the coefficient of $x^0$ in $x \frac{\partial}{\partial x} f(x)$ is zero. So:
\begin{equation}
	\begin{array}{l}
	0 = [x^0]\ x \frac{\partial}{\partial x} \frac{w(x)^{m+1}}{x^k} \\
	= [x^0]\ \left ( -k(a_0 + a_1 x + \ldots + a_N x^N) \frac{w(x)^m}{x^k} \right . \\ \hspace{2em} \left . + (m+1)(a_1 + 2 a_2 x + \ldots + N a_N x^{N-1} \frac{w(x)^m}{x^{k-1}} \right ) \\
	= [x^0]\ \left ( -k \left [ a_0 \frac{w(x)^m}{x^k} + a_1 \frac{w(x)^m}{x^{k-1}} + \ldots + a_N \frac{w(x)^m}{x^{k-N}} \right ] \right . \\ \left . \hspace{2em} + (m+1) \left [ a_1 \frac{w(x)^m}{x^{k-1}} + 2 a_2 \frac{w(x)^m}{x^{k-2}} + \ldots + N a_N \frac{w(x)^m}{x^{k-N}} \right ] \right ) \\
	= -k [ a_0 c_k + a_1 c_{k-1} + \ldots + a_N c_{k-N} ] \\ \hspace{2em} + (m+1) [ a_1 c_{k-1} + 2 a_2 c_{k-2} + \ldots + N a_N c_{k-N} ] \\
	\end{array}
\end{equation}
From which the recurrence follows. The standard notation~\cite{Wilf1993} that $[x^n] f(x)$ is the coefficient of $x^n$ in $f(x)$ has been used.

\section{Exponentiation in the DTM}
Combining the recurrence relation with the DTM formalism is straightforward. Note that $W(k)$ is the the $k^{\hbox{th}}$ Taylor-series coefficient of $w(x)$, and this yields the table entry in Table~\ref{tab:dtmexp}. Care must be taken, however, to insure that if $Y(0)$ is zero then the replacement $y(x) \rightarrow x \bar{y}(x)$ where $\bar{y}(x) = \frac{y(x)}{x}$ is made.

\begin{table}[h]
\begin{center}
\begin{tabular}{|p{4cm}|p{6cm}|}
\hline
\bf{Original Function}                  & \bf{Transformed Function}                                                           \\
\hline
$w(x) = y(x)^n$, $Y(0) \neq 0$ & $W(k) = \frac{1}{kY(0)} \sum_{j=1}^k [(m+1)j - k] Y(j) W(k-j)$ \newline $W(0) = Y(0)^m$ \\
$w(x) = e^{y(x)}$              & $W(k) = \frac{1}{k} \sum_{j=1}^k j Y(j) W(k-j)$ \newline $W(0) = e^{Y(0)}$ \\
\hline
\end{tabular}
\end{center}
\caption{Exponentiation operations under DTM. $n$ is a positive integer.}
\label{tab:dtmexp}
\end{table}

Using this formulation greatly simplifies the expressions resulting from exponentiating a function compared to the previously-available method~\cite{Hassan2007}. Specifically, for $w(x) = y(x)^m$ the representation constructed iteratively from the multiplication rule was
\begin{equation}
	\label{eq:oldway}
	W(k) = \sum_{k=0}^\infty \sum_{k_{m-1}=0}^k \sum_{k_{m-2}=0}^{k_{m-1}} \cdots \sum_{k_2=0}^{k_3} \sum_{k_1=0}^{k_2} Y(k_1) Y(k_2 - k_1) \cdots Y(k_{m-1} - k_{m-2}) Y(k - k_{m-1}) \, .
\end{equation}
Numerically, there is a significant advantage to using the recurrence in Table~\ref{tab:dtmexp}, which requires $\mathrm{O}(N)$ computations, compared to the iterative method (Equation~\ref{eq:oldway}), which requires $\mathrm{O}(m N \log N)$ computations.

\section{The exponential function in DTM}
Using Zeilberger's method, it is possible to derive a recurrence for $e^{w(x)}$:
\begin{equation}
	\begin{array}{l}
	0 = [x^0]\ x \frac{\partial}{\partial x} \frac{e^{w(x)}}{x^k} \\
	= [x^0]\ \left ( -k \frac{e^{w(x)}}{x^k} + (a_1 + 2 a_2 x + \ldots + N a_N x^{N-1} \frac{e^{w(x)}}{x^{k-1}} \right ) \\
	= [x^0]\ \left ( -k \frac{e^{w(x)}}{x^k} + \left [ a_1 \frac{e^{w(x)}}{x^{k-1}} + 2 a_1 \frac{e^{w(x)}}{x^{k-2}} + \ldots + N a_N \frac{e^{w(x)}}{x^{k-N}} \right ] \right ) \\
	= -k c_k + [ a_1 c_{k-1} + 2 a_2 c_{k-2} + \ldots + N a_N c_{k-N} ] \\
	\end{array}
\end{equation}
The resulting recurrence is also listed in Table~\ref{tab:dtmexp}. Again, it is much simpler than the iterative exponentiation construction in combination with the Taylor-series expansion of the exponential function~\cite{Hassan2007}
\begin{equation}
	\label{eq:expoldway}
	W(k) = \sum_{m=0}^\infty \frac{1}{m!} \sum_{k=0}^\infty \sum_{k_{m-1}=0}^k \sum_{k_{m-2}=0}^{k_{m-1}} \cdots \sum_{k_2=0}^{k_3} \sum_{k_1=0}^{k_2} Y(k_1) Y(k_2 - k_1) \cdots Y(k_{m-1} - k_{m-2}) Y(k - k_{m-1}) \, .
\end{equation}

\section{An example: The One-Dimensional Planar Bratu Problem}
To exemplify the formulation presented here, the simplified recurrence will be applied to the one-dimensional planar Bratu problem which had previously been studied using DTM by Hassan and Ert\"urk~\cite{Hassan2007}. The differential equation is
\begin{equation}
\begin{array}{ll}
u'' + \lambda e^u = 0 & 0 \leq x \leq 1 \\
u(0) = 0, u(1) = 0 & \\
\end{array}
\end{equation}
and the previously-derived DTM solution is
\begin{equation}
	\label{eq:oldsltn}
	\begin{array}{l}
	U(k+2) = -\frac{-\lambda}{(k+1)(k+2)} \times \\
	\hspace{2em} \left ( \sum_{m=0}^\infty \frac{1}{m!} \sum_{k=0}^\infty \sum_{k_{m-1}=0}^k \sum_{k_{m-2}=0}^{k_{m-1}} \cdots \sum_{k_2=0}^{k_3} \sum_{k_1=0}^{k_2} \right . \\ \hspace{4em} \left . U(k_1) U(k_2 - k_1) \cdots U(k_{m-1} - k_{m-2}) U(k - k_{m-1}) \right )
	\end{array}
\end{equation}
for which the coefficients $\{U(k)\}$ are computed starting from
\begin{equation}
	\begin{array}{l}
	U(0) = 0 \\
	U(1) = \gamma \\
	\end{array}
\end{equation}
where $\gamma$ is a constant fixed by the $x=1$ boundary condition
\begin{equation}
	\sum_{k=0}^\infty U(k) = 0 \, .
\end{equation}
This can be compared to the analytic solution~\cite{Hassan2007}
\begin{equation}
	u(x) = -2 \ln \left [ \frac{\cosh \left ( \left ( x - \frac{1}{2} \right ) \frac{\theta}{2} \right )}{\cosh \left ( \frac{\theta}{4} \right ) } \right ]
\end{equation}
where $\theta$ solves
\begin{equation}
	\label{eq:thcond}
	\theta = \sqrt{2 \lambda} \cosh \left ( \frac{\theta}{4} \right ) \, .
\end{equation}
Depending on the value of $\lambda$, Equation~\ref{eq:thcond} has zero, one or two solution(s). See the above-referenced paper by Hassan and Ert\"urk, and references cited therein, for further details.

Applying the recurrence for $e^{u(x)}$ provided here yields the following DTM solution
\begin{equation}
	\begin{array}{l}
	U(k+2) = -\frac{-\lambda}{(k+1)(k+2)} W(k) \\
	W(k) = \frac{1}{k} \sum_{j=1}^k j U(j) W(k-j) \\
	W(0) = e^{U(0)}
	\end{array}
\end{equation}
which can be simplified by writing $W(k) = -\frac{(k+1)(k+2)}{\lambda} U(k+2)$ to
\begin{equation}
	\label{eq:newsltn}
	\begin{array}{ll}
	U(k+2) = \frac{1}{k(k+1)(k+2)} \sum_{j=1}^k j (k-j+1)(k-j+2) U(j) U(k-j+2) & k \geq 1 \\
	U(2) = -\frac{\lambda}{2} e^{U(0)} & \\
	\end{array}
\end{equation}
The first few values of the recurrence are
\begin{equation}
	\begin{array}{l}
	U(0) = 0 \\
	U(1) = \gamma \\
	U(2) = -\frac{\lambda}{2} \\
	U(3) = -\gamma \frac{\lambda}{6}
	\end{array}
\end{equation}
and these are the same regardless of whether Equation~\ref{eq:oldsltn} or Equation~\ref{eq:newsltn} is used.

\section{Conclusion}
The rule derived from J.C.P. Miller's recurrence, and a similar recurrence derived for the exponential function, given in Table~\ref{tab:dtmexp}, are far simpler than the representations previously used in the literature (Equations~\ref{eq:oldway} and~\ref{eq:expoldway}). It should prove easier to apply DTM, and any power-series-solution technique, to nonlinear differential equations using the recurrence relations presented here.

\section{Acknowledgments}
I thank the United States Department of Energy Computational Science Graduate Fellowship, provided under grant DE-FG02-97ER25308, for supporting this research. I also thank Richard Easther for providing valuable feedback and suggestions.

\bibliographystyle{plain}
\bibliography{dtm_miller}

\end{document}